\documentstyle[12pt]{article}
\newcommand{\h}{ \cal H}
\newcommand{\fti}{\{f_i \}_{i \in I}}
\newcommand{\ftj}{ \{f_i\}_{i \in J}}

\newcommand{ \sijpc}{\sum_{ i \in I- \jp}}
\newcommand{\sjn}{ \sum_{j=1}^n}

\newcommand{\sui}{\sum_{i \in I}}
\newcommand{\suj}{\sum_{i \in J}}
\newcommand{\jp}{ J^{\prime}}

\newcommand{\fc}{ <f,f_i >}
\newcommand{\fct}{ | <f,f_i> |^2}

\newcommand{\nft}{ ||f||^2 }
\newcommand{\si}{S^{-1}}
\newcommand{\fei}{ <f,e_i >}
\newcommand{\feo}{ <f, e_1>}

\newcommand{\ytu}{\{y_i \}_{i=1}^{\infty}}
\newcommand{\suu}{\sum_{i=1}^{\infty}}
\newcommand{\ei}{{\epsilon}_i}
\newcommand{\eij}{{\ei}^j}
\newcommand{\snj}{ \sum_{i=n_j}^{n_{j+1}-1}}

\newcommand{\ftu}{\{f_i \}_{i=1}^{\infty}}
\begin{document}
\title{ Hilbert space frames containing a Riesz basis and Banach spaces
which have no subspace isomorphic to $c_0$.}
\author{Peter G. Casazza and Ole Christensen 
\\ This paper is dedicated to Ky Fan.
\thanks{ The first author
acknowledges support by NSF Grant DMS-9201357 and a grant from the Danish
Research Foundation.  1993 Mathematics Subject Classification:
Primary 42C99, 46B03, 46C99. }}
\maketitle
\begin{abstract} We prove that a Hilbert space frame $\fti$ contains a Riesz 
basis if every subfamily $\ftj ,
J \subseteq I ,$ is a frame for its closed span. Secondly we give
a new characterization of Banach spaces 
which do not have any subspace isomorphic to $c_0$. This result
immediately leads to an improvement of a recent theorem of
Holub concerning frames consisting of a Riesz basis plus finitely many 
elements.  
\end{abstract}
\section{Introduction.}
Let $\h$ be a separable Hilbert space with the inner product
$< \cdot , \cdot > $ linear in the first entry. A family $\fti$ of
elements in $\h$ is called a {\it frame} for $\h$ if
$$\exists A,B>0 : \ \ A \nft \leq \sui \fct \leq B \nft , \ \
\forall f \in \h .$$
\ $A$ and $B$ are called {\it frame bounds}. If $\fti$ is a frame, then
\\ $Sf:= \sui <f,f_i >f_i $ defines a bounded invertible operator on $\h$.
This fact leads to the {\it frame decomposition}: every $f \in \h$ has
the representation 
$$ f= SS^{-1}f= \sui <f, \si f_i >f_i .$$
So a frame can be considered as a generalized basis in the sense
that every element in $\h$ can be written as a linear combination 
of the frame elements. Actually one has that
$$\fti \ \mbox{ is a Riesz basis} \Leftrightarrow [ \fti \ \mbox{is 
a frame and} \ \sui c_i f_i =0 \Rightarrow c_i =0, \  \forall i .]$$
Recently more authors have been interested in the relation between
frames and Riesz bases. Holub \lbrack H\rbrack \ concentrates on
{\it near-Riesz bases}, i.e., frames consisting of a Riesz basis plus 
finitely many elements. He is able to give equivalent 
characterisations of such frames. Seip \lbrack Se\rbrack \ deals only
with frames of complex exponentials. Among his very interesting
results one finds examples of frames which do not contain a
Riesz basis. On the other hand he proves that all frames which have 
appeared ``naturally'' until his paper contains a Riesz basis. \\
Using different techniques, the present authors have constructed
a frame not containing a Riesz basis \lbrack CC\rbrack . Furthermore
one of the authors gave the first condition implying that a frame
contains a Riesz basis \lbrack C1\rbrack . The purpose of the first
part of the present paper is to show that the conclusion is true
under weaker conditions than in \lbrack C1\rbrack . We also give an
example where the new theorem can be used, but where the condition
in the version from \lbrack C1\rbrack \ is not satisfied. \\ \\

In the second part of the paper we present a new characterization of
Banach spaces which do not have
any subspace isomorphic to $c_0$. This is important in itself, but the
reason for combining it with the frame result above is that it immediately
leads to an improvement of a recent frame result of Holub  
\lbrack H\rbrack. Holub shows that if a frame is norm-bounded below,
then it is a near-Riesz basis if and only if it is {\it unconditional},
which means that if a series $\sui c_i f_i $ converges, then it converges
unconditionally. Actually we are able to prove  the same without
the condition about norm-boundedness; however, this property follows as a
consequence of the situation. 
\section{Frames containing a Riesz basis.}
A frame $\fti$ is called a {\it Riesz frame} if every subfamily
$\ftj$ is a frame for its closed span, with a lower bound $A$ common
for all those frames. One of the main results in \lbrack C1\rbrack \ is
\\ \\ {\bf Theorem 2.1:} {\it Every Riesz frame contains a Riesz basis.} \\ \\
The main ingredient in the proof is an application of Zorn's lemma.
Our goal here is to show that the conclusion actually holds 
without the assumption about a common lower bound. This is important 
in practice, since one might be in the situation that some theoretical
arguments give the frame property, but no knowledge of the bounds. 
However the proof of this more general theorem
is much more complicated and in fact Theorem 2.1 is part of the 
results we use in the proof. We need a lemma: \\ \\
{\bf Lemma 2.2:} {\it Let $\fti$  be a frame. Given $\epsilon > 0 $ and a
finite set $J \subseteq I $, there exists a finite set $ \jp $ containing
$J$ such that  
$$ \sijpc | \fc |^2 \leq \epsilon \cdot \nft , \ \forall f \in span \ftj .$$}
{\bf Proof:} Let $\{e_j \}_{j=1}^n $ be an orthonormal basis for $span \ftj .$
Given $ \epsilon >0 $ we take the index set $\jp $ containing $I$ such that
$$ \sijpc | <e_j , f_i > |^2 \leq \frac{\epsilon}{n}, \ 
\mbox{for all} \ j \in J .$$
Now take $f \in span \ftj .$ Writing $f = \sum_{j=1}^n <f,e_j >e_j $
we get
$$ \sijpc | \fc |^2 =  \sijpc | \sjn <f,e_j >< e_j , f_i > |^2 $$ \ $$\leq 
\sijpc \sjn | <f,e_j> |^2 \sjn |<e_j , f_i >|^2 
$$ \ $$ = \sijpc \sjn |< e_j, f_i >|^2 \cdot \nft \leq 
\epsilon \cdot \nft .$$
{\bf Proposition 2.3:}
{\it Let $\{f_i \}_{i \in I}$ be a frame with the property that every subset of
$\{f_i \}$ is a frame for its closed linear span.  Then there is
an $\epsilon > 0,$ and finite subsets $J \subset J' \subset I,$
with the property:  For every $J" \subset I-J',$ the family
$\{f_i \}_{i \in J \cup J"}$ has lower frame bound $\ge \epsilon.$}
\\ \\ {\bf Proof:}
We assume the proposition fails and construct by induction  
sequences of finite subsets $J_1,J_2,\ldots$ and  
$J_{1}^{'},J_{2}^{'},\ldots$
with the following properties.

(1)  ${\cup}_{j=1}^{n} J_j \subset J_{n}^{'},$

(2)  For every $f \in \ {\mbox{span}} \{f_i \}_{i \in {\cup}_{j=1}^{n} J_j},$
with $\|f\| = 1$ and
$$
{\sum_{i \in I-J_{n}^{'}} |<f,f_i>|^{2}} \le \frac{1}{n},
$$

(3)  There is some $f \in \ {\mbox{span}} \{f_i \}_{i \in {\cup}_{j=1}^{n}  
J_j},$ with $\|f\| = 1$ and
$$
{\sum_{i \in {\cup}_{j=1}^{n}J_j} |<f,f_i>|^2} \le \frac{1}{n}.
$$

We will quickly check the induction step.  Assume  
$J_1,J_2,\ldots,J_n$ and $J_{1}^{'},J_{2}^{'},\ldots,J_{n-1}^{'}$  
have been chosen to satisfy (1)-(3) above.  By lemma 2.2, there is a  
finite set $J_{n}^{'}
\subset I$ with ${\cup_{j=1}^{n} J_j} \subset J_{n}^{'}$ satisfying  
(2) above with the constant 1/(n+1).  Given $\epsilon =$ 1/(n+1),  
and
$J = \cup_{j=1}^{n} J_j$ and $J \subset J_{n}^{'},$ our assumption  
that the proposition fails implies there is a finite set $J_{n+1}  
\subset I-J_{n}^{'}$ so that (3) holds for 1/(n+1).  This completes  
the induction.  We now let $J = \cup_{n=1}^{\infty} J_n.$  It is  
easily seen from (2) and (3) above that $\{f_i \}_{i \in J}$ is not a  
frame for its closed linear span.  That is, for each n, there is a  
$f \in \mbox{span} \{f_i \}_{i \in \cup_{j=1}^{n} J_j} $ with $\|f\| = 1$  
and satisfying (3). We now have by (2),   
$$
\sum_{i \in J} |<f,f_i>|^2 =
$$
$$
{\sum_{i \in {\cup_{j=1}^{n}}J_j}|<f,f_i>|^2} + {\sum_{i \in  
{\cup_{j=n+1}^{\infty}}J_j}|<f,f_i>|^2} \le \frac{1}{n} +  
\frac{1}{n} = \frac{2}{n}. $$
This contradiction completes the proof of the proposition. \\ \\
{\bf Theorem 2.4:} {\it If every subset of $\fti$ is a frame for its closed
linear span, then $\fti$ contains a Riesz basis.} \\ \\
{\bf Proof:} 
By proposition 2.3, there exists an $\epsilon > 0$ and a finite sets
J and J' with $J \subset J'$ so that whenever $J'' \subset I-J'$  
the lower frame bound of $\{f_i \}_{i \in J \cup J''}$ is $\ge  
\epsilon.$
Let $P$ denote the orthogonal projection of $H$ onto \ $\mbox  
{span}\{f_i \}_{i \in J} $  Then for all $J'' \subset I-J',$ if $f \in  
\mbox{span} \{ (I-P)f_i \}_{i \in J''}$ ,then $f \in 
\mbox{span} \{f_i \}_{i \in J  \cup J'} ,$ so
$$ \epsilon ||f||^2  \leq \sum_{i \in J \cup J''} |<f,f_i>|^2 
= \sum_{i \in J \cup J''} |<(I-P)f,f_i>|^2 $$ \
$$ = \sum_{i \in J''} |<(I-P)f,f_i>|^2 \\
= \sum_{i \in J''} |<f,(I-P)f_i>|^2. $$

It follows that for every $J'' \subset I-J',$ the set  
$\{(I-P)f_i \}_{i \in J''}$ has lower frame bound $\epsilon > 0.$
Obviously every frame 
$\{(I-P)f_i \}_{i \in \jp}$ has the same upper bound as $\fti$, so Theorem 2.1
applied to $\{(I-P)f_i \}_{i \in I-J^{\prime}}$ 
shows that there exists a subset $I^{\prime} \subseteq I-J^{\prime}$ such that 
$\{(I-P)f_i \}_{i \in I^{\prime}}$ is a Riesz basis for its closed span . 
Hence, for all sequences $\{c_i \} \in l^2(I^{\prime})$,
$$ ||\sum_{i \in I^{\prime}} c_i f_i || \geq 
|| \sum_{i \in I^{\prime}} c_i (I-P)f_i || \geq \epsilon 
\sqrt{ \sum_{i \in I^{\prime}} |c_i |^2 },$$ i.e., 
$\{f_i \}_{i \in I^{\prime}} $  is a Riesz basis for its closed span.
But since $dim(P \h) < \infty$ it can be extended to a Riesz basis for $\h$
just by adding finitely many elements. \\ \\   
To prove that Theorem 2.4 really is an improvement of Theorem 2.1 one needs
an example of a frame, where every subfamily is a frame for its closed
span, but where there does not exists a common lower bound for all those
frames. We present such an example now: \\ \\
{\bf Example:} Let $\{e_i \}_{i=1}^{\infty}$ be an orthornormal basis
for $\h$ and define
$$\{ f_i \}_{i \in K} := \{e_i, e_i + \frac{1}{2^i}e_1 \}_{i=2}^{\infty}.$$
First we show that every subfamily $\{f_i \}_{i \in L}$ is a frame for
its closed span. For convenience, write the subfamily as
$$ \{f_i \}_{i \in L} = \{e_i \}_{i \in I} \cup 
\{e_i + \frac{1}{2^i}e_1 \}_{i \in J} .$$
First we assume that $I \cap J = \emptyset .$ Then 
$\{e_i \}_{i \in I \cup J} $ is an orthonormal basis for its closed
span. The idea is now to show that $\{f_i \}_{i \in L}$ is a perturbation
of this family and thereby conclude that the family itself is a frame. Since 
$$ \sum_{i \in I \cup J} ||f_i -e_i ||^2 = \suj [\frac{1}{2^i}]^2
\leq \sum_{i=2}^{\infty} [\frac{1}{2^i}]^2 < 1$$
we conclude by \lbrack C1, Cor. 2.3. b)\rbrack \ that $\{f_i \}_{i \in L}$
is  a Riesz basis for its closed span, as desired.
\\ \\ Now assume that $I \cap J \neq \emptyset .$ Write
$$ \{f_i \}_{i \in L}= \{e_i, e_i + \frac{1}{2^i}e_1 \}_{i \in I \cap J}
\cup \{e_i \}_{i \in I-J} \cup \{e_i + \frac{1}{2^i}e_1 \}_{i \in J-I} .$$
Clearly $\overline{span} \{f_i \}_{i \in L}=
\overline{span} \{ \{e_1 \} \cup \{e_i \}_{i \in I \cup J} \}$. The sequence 
$ \{f_i \}_{i \in L}$ is a Bessel sequence (i.e., the upper frame condition 
is satisfied) so by \lbrack C2, Cor. 4.3\rbrack \ we are done if we can
show that the operator 
$$ T : l^2(L) \rightarrow \overline{span} \{f_i \}_{i \in L},
 \ \ T \{c_i \} = \sum_{i \in L} c_i f_i $$ is surjective. Now, let $f \in 
\overline{span} \{f_i \}_{i \in L} .$ We want to write $f$ as a linear 
combination of elements $f_i$ with coefficients from $l^2(L)$. First,  
$$ f = \sum_{i \in I-J} \fei e_i + \sum_{i \in J-I} \fei e_i +
\sum_{i \in I \cap J} \fei e_i + \feo e_1 $$ \ $$ 
= \sum_{i \in I-J}\fei e_i + \sum_{i \in J-I} \fei (e_i + \frac{1}{2^i}e_1)
$$ \ $$ - \sum_{i \in J-I} \fei \frac{1}{2^i}e_1 + \feo e_1 +
\sum_{i \in I \cap J} \fei e_i .$$
Choose $ n \in I \cap J .$ Then 
$$f = \sum_{i \in I -J} \fei e_i  + 
\sum_{i \in J-I} \fei (e_i + \frac{1}{2^i}e_1) $$\ $$+ 
(-2^n \sum_{i \in J-I} \fei \frac{1}{2^i} + 2^n \feo)(e_n + \frac{1}{2^n}e_1) 
$$ \ $$-(-2^n \sum_{i \in J-I} \fei \frac{1}{2^i}+ 2^n \feo)e_n $$ \ $$
+ <f, e_n> e_n + \sum_{i \in I \cap J- \{n\}} \fei e_i  $$ \ $$
= \sum_{i \in I \cup J} \fei e_i + 
\sum_{i \in J-I} \fei (e_i +\frac{1}{2^i}e_1) $$ \ $$
+ (-2^n \sum_{i \in J-I} \fei \frac{1}{2^i} + 2^n \feo ) \cdot
(e_n + \frac{1}{2^n}e_1) $$ \ $$ + (2^n \sum_{i \in J-I} \fei \frac{1}{2^i}
- 2^n \feo + <f, e_n>)e_n + \sum_{i \in I\cap J - \{n\}} \fei e_i .$$
So every $f \in \overline{span} \{f_i \}_{i \in L}$ can be written as a linear 
combination of the elements in $\{f_i \}_{i \in L}$ and
obviously the coefficient sequence is in $l^2(L)$. That is, $T$ is surjective,
and the proof is complete. Now we show that there is no common lower bound
for all subframes. Let  $n \in N$ and 
consider the family $\{ e_n , e_n + \frac{1}{2^n}e_1 \}$,
which is a frame for $span \{e_1, e_n \} .$ Since 
$$ | <e_1 , e_n >|^2 + |< e_1, e_n + \frac{1}{2^n}e_1 >|^2 = 
\frac{1}{2^{2n}} \cdot || e_1 ||^2 $$ the lower bound for this frame is at most
$\frac{1}{2^{2n}}$. Hence there is no common lower bound. 
\section{Banach spaces isomorphic to $c_0$ and near-Riesz bases
in Hilbert spaces.}
In this section we prove the following: \\ \\
{\bf Theorem 3.1:} {\it Let $\ftu$ be a frame for the Hilbert space $\h$.
Then $\ftu$ is unconditional if and only if it is a near-Riesz basis.} \\ \\
With the additional assumption that the $f_i's$ are norm bounded below the
result is proven by Holub \lbrack H\rbrack \ using a result of 
Heil \lbrack He, p.168\rbrack . Theorem 3.1 shows that this assumption is
superfluous. However, it is a consequence of the situation, since every
near-Riesz basis is norm-bounded below. 
The ``if ' above follows from this property and the original result of
Holub. The ``only if'' part is more complicated, but 
actually we prove a much more general result
concerning abstract Banach spaces. This result has independent interest
in Banach space theory in that it classifies those Banach spaces which
do not have any subspace isomorphic to $c_o $. \\ \\
{\bf Theorem 3.2:} {\it Let $X$ be any Banach space.
The following are equivalent: \\ \\
(1) No subspace of $X$ is isomorphic to $c_0$. \\ \\
(2) If $\ytu \subseteq X$ is a sequence so that whenever 
$ \suu a_i y_i $ converges for some coefficient sequence $\{a_i \} $, the
series must converges unconditionally, then there is some $n_o \in N$
so that $\{y_i \}_{i=n_0}^{\infty}$ is a unconditional basis for its closed
span.} \\ \\
The ``only if'' part of Theorem 3.1 is a consequence of Theorem 3.2. 
Actually a Hilbert space satisfies (1) of Theorem 3.2, so if 
$\ftu$ is a near-Riesz basis for a Hilbert space $\h$, 
then there exists a $n_0$ such that
$\{f_i \}_{i=n_0}^{\infty}$ is a unconditional basis for its closed span
and a frame for its closed span.
Here we used the fact that if one deletes finitely many elements from a 
frame, then one still has a family which is a frame for its closed span
 (\lbrack C3, Lemma 2\rbrack \ or \lbrack CH, Theorem 7\rbrack \ for a more 
general statement). So by the characterization of Riesz bases in the 
introduction, $\{f_i \}_{i=n_0}^{\infty}$ is a Riesz basis for its closed 
span. This space has finite codimension, so adding finitely 
many elements we obtain a Riesz basis for $\h$, and Theorem 3.1 follows. \\ \\
The proof of Theorem 3.2 requires some preliminary results. \\ \\
{\bf Lemma 3.3:} {\it If $X$ is a Banach space, $\ytu \subseteq X$ and
$\{c_i \}_{i=1}^{\infty} $ are scalars so that $ \suu c_i y_i $
converges unconditionally in $X$, then
$$ lim_{n \rightarrow \infty} sup_{\ei = \pm 1} || \sum_{i=n}^{\infty} 
\ei c_i y_i || =0 .$$} 
{\bf Proof:} If the conclusion of the lemma fails, then there is some 
$\epsilon >0$ and natural numbers $n_1 < n_2 < ...$ and some $\eij = \pm 1 ,
j=1,2,... \ \mbox{and} \
i= n_j , n_j+1...,$ with $ || \sum_{i=n_j}^{\infty} \eij c_i y_i || \geq
\epsilon . $ Since $\suu c_i y_i $ converges unconditionally,
$ \sum_{i=n_j}^{\infty} \eij c_i y_i $ converges in $X$, and hence
$ lim_{k \rightarrow \infty} || \sum_{i=k}^{\infty} \eij c_i y_i || =0 $
for all $j$. Therefore, by switching to a subsequence of $n_j$ 
(let us call it $n_j$ again), we have for all $j=1,2....,$ \
$$ || \snj \eij c_i y_i || \geq \frac{\epsilon}{2} .$$
Letting $d_i := \eij c_i, \ \mbox{for} \ \ n_j  \leq i \leq n_{j+1}-1 $
we have that 
$$ \sum_{j=1}^{\infty} \snj \eij c_i y_i = \sum_{i=n_1}^{\infty} d_i y_i $$
converges in $X$ , since it is just a change of signs on the unconditionally
convergent series  $\sum_{i=n_1}^{\infty} c_i y_i .$ However, 
$$ || \snj d_i y_i || = || \snj \eij c_i y_i || \geq \frac{\epsilon}{2} $$
implies that $\sum_{i=n_1}^{\infty} d_i y_i$ does not converges in $X$. 
This contradiction completes the proof of Lemma 3.3. \\ \\
Next we introduce the notation needed for the proof of Theorem 3.2. If
$\{x_i \}_{i=1}^{\infty}$ is a basis for a Banach space 
$X, \ n_1 < n_2 < ...$ are natural numbers, and $y_j = \snj c_i x_i $
are vectors in $X$, we call $\ytu$ a {\it block basic sequence} of
$\{x_i \}_{i=1}^{\infty}$. If $\{x_i \}_{i=1}^{\infty}$ is a unconditional
basis for $X$, then it is easily seen that a block basic sequence 
$\ytu$ is a unconditional
basis for its closed span with unconditional basis constant less than or equal
to the unconditional basis constant for $\{x_i \}_{i=1}^{\infty} $ in $X$. 
A series $\sum_{n} x_n$ in a Banach space is said to be {\it weakly  
unconditionally Cauchy} if given any permutation $\pi$ of the  
natural numbers,
$(\sum_{k=1}^{n} x_{\pi (k)})$ is a weakly Cauchy sequence;  
alternatively,
$\sum_{n} x_n$ is weakly unconditionally Cauchy if and only if for each  
$x^* \in X^*$, $\sum_{n}|x^{*}(x_n)| < \infty$ (see \lbrack D \rbrack,  
Chapter V).  Let $\hat N$ denote the family of all finite subsets of  
the natural numbers.  We will need a result which follows  
immediately from theorem 6, page 44, of \lbrack D\rbrack . \\ \\
{\bf Proposition 3.4}  {\it The following statements are equivalent: \\
(1)  $\sum_{n} x_n$ is weakly unconditionally Cauchy.
\\ \\
(2) $ \mbox{sup}_{\Delta \in \hat N}|| \sum_{n \in \Delta} x_n || < \infty. $}
\\ \\Also, we will need Theorem 8, page 45, from \lbrack D\rbrack, which 
we now state for completeness. \\ \\
{\bf Proposition 3.5}
{\it Let $X$ be a Banach space.  Then, in order that each weakly  
unconditionally
Cauchy series in $X$ be unconditionally convergent, it is necessary  
and sufficient that $X$ contains no copy of $c_0$.} \\ \\ 
Now we are ready to prove Theorem 3.2.  \\
$(2)  \Rightarrow (1)$: It suffices to show that $c_o$ fails property (2).
Let $ \{e_n \}_{n=1}^{\infty}$ be the unit vector basis of $c_o$ and define
$$ y_{2n} = e_n, \ \ y_{2n+1} = e_n , \ n=1,2,.. $$
We will show that $\{y_n \}$ satisfies the hypotheses of (2) but fails the
conclusion. So assume that $ \sum_{n=1}^{\infty} c_n y_n $ converges in 
$c_0 $. Since $ ||y_n|| =1 $, for all $n$, it follows that 
$lim_{n \rightarrow \infty} |c_n | =0 $. Given any ${\epsilon}_n = \pm 1 , \
 || \sum_{n=m}^{\infty} {\epsilon}_n c_n y_n ||_{\infty} \leq 2 sup_{m \leq k}
| {\epsilon}_k c_k | $. Hence, $$ lim_{ m \rightarrow \infty} || 
\sum_{n=m}^{\infty} {\epsilon}_n c_n y_n ||_{\infty}
\leq 2 lim_{m \rightarrow \infty} sup_{m \leq k} |c_k | =0 .$$
So \ $\sum_{n=1}^{\infty} {\epsilon}_n c_n y_n $ converges in $c_o $.
That is, whenever $\sum_{n=1}^{\infty} c_n y_n$ converges in $c_0$, then the 
series converges unconditionally. So the hypotheses of (2) in Theorem 3.2
are satisfied. But clearly the conclusion of (2) fails since any subset of
$\{y_n \}$, which contains all but a finite number of the $y_n$, must 
contain two equal elements and hence cannot be independent. \\ \\

$(1) \Rightarrow (2).$  We proceed by way of contradiction.  So  
assume (1) and the hypotheses of (2) are satisfied, but the  
conclusion of (2) fails.
Alternately applying this assumption and lemma 3.3, we find natural  
numbers $n_1, n_2, \ldots$, and ${\epsilon}_{i}^{j} = \pm 1$, and  
scalars $\{c_{i}\}_{i = 1}^{\infty}$ so that for all j,
\vskip10pt
(3)  $\|\sum_{i=n_{j}+1}^{n_{j+1}} c_iy_i\| < \frac{1}{2^j},$
\vskip10pt
(4)  $\frac{1}{2} \le \frac{1}{2}{\mbox{sup}}_{{\epsilon}_i = \pm  
1}\|\sum_{i=n_{j}+1}^{n_{j+1}} {\epsilon}_{i} c_{i} y_{i}\|
\le \|\sum_{i=n_{j}+1}^{n_{j+1}} {\epsilon}_{i}^{j} c_{i} y_{i}\| = 1.$ 
\vskip10pt

We let $z_j = \sum_{i=n_{j}+1}^{n_{j+1}} {\epsilon}_{i}^{j} c_i y_i$, for
$j = 1, 2, \ldots$. \\ 

{\bf Claim:}  ${\mbox{sup}}_{\Delta \in {\hat N}} \|\sum_{i \in  
\Delta} z_i\| = \infty$.
\\ \\
This claim follows quickly.  That is, if this sup was finite, then  
$\sum_{i} z_i$ would be unconditionally Cauchy by Proposition 3.4.  
But
 then since we assumed $c_0$ does not embed into our space, by  
Proposition 3.5, it would follow that this series is unconditionally  
convergent.  But this is ridiculous since $\| z_i \| = 1$, for all  
$i = 1, 2, \ldots$.  This completes the proof of the Claim.

  By applying our claim and choosing successive subsets $\Delta \in  
{\hat N}$, and reindexing we have the following:  There are natural  
numbers $n_1, n_2, \ldots$, natural numbers $0 = m_{0} < m_1 < m_2  
< \cdots$, a sequence of scalars $\{c_i\}$ and choices of signs  
${\epsilon}_{i}^{j} = \pm 1$, so that
\vskip10pt
(5)  $\|\sum_{i=n_{j}+1}^{n_{j+1}} c_{i}y_{i}\| < \frac{1}{2^j},$
\vskip10pt
(6)  $\frac{1}{2} \leq \frac{1}{2}{\mbox{sup}}_{{\epsilon}_i = \pm  
1}\|\sum_{i=n_{j}+1}^{n_{j+1}} {\epsilon}_{i} c_{i} y_{i}\|
\leq \|\sum_{i=n_{j}+1}^{n_{j+1}} {\epsilon}_{i}^{j} c_{i} y_{i}\| = 1.$
\vskip10pt
(7)  $\|\sum_{j=m_{k}+1}^{m_{k+1}}(\sum_{i=n_{j}+1}^{n_{j+1}}  
{\epsilon}_{i}^{j} c_{i} y_{i}) \| = K_{k} \geq k.$
\vskip10pt

We will now show that the series
$$ (8) \ \ \ \ 
\sum_{k=1}^{\infty}  
{\frac{1}{K_k}}[\sum_{j=m_{k}+1}^{m_{k+1}}(\sum_{n_{j}+1}^{n_{j+1}}  
c_{i}y_{i})] $$
converges in $X$ as a series in $c_{i}y_{i}$, but the series does  
not converge unconditionally.

That the series does not converge unconditionally can be proven easily.
For any $k = 1,2,\ldots$,
$$
{\mbox{sup}}_{{\epsilon}_{i} = \pm 1}  
{\frac{1}{k_k}}\|\sum_{j=m_{k}+1}^{m_{k+1}}(\sum_{i=n_{j}+1}^{n_{j+1}}  
{\epsilon}_{i} c_{i}y_{i})\| \ge  
{\frac{1}{K_k}}\|\sum_{j=m_{k}+1}^{m_{k+1}}(\sum_{i=n_{j}+1}^{n_{j+1}}  
{\epsilon}_{i}^{j} c_{i}y_{i})\| = 1. $$
That is, the series in (8) fails lemma 3.3 and hence is not unconditionally
convergent.
To prove that the series in (8) converges, we must check that the  
"tail end" of the series converges to $0$ in norm.  So consider  
$\sum_{i=s}^{\infty} c_{i}y_{i}$, and fix k with $m_{k}+1 \leq l \leq  
m_{k+1}$ and $n_{l}+1 \leq s \leq n_{l+1}$.  Then
$$
\|{\frac{1}{K_t}}{\sum_{i=s}^{n_{l+1}}c_{i}y_{i}} +  
{\frac{1}{K_k}}\sum_{j=l+1}^{m_{k+1}}(\sum_{i=n_{j}+1}^{n_{j+1}}  
c_{i}y_{i}) +  
\sum_{t=k+1}^{\infty}\frac{1}{K_t}[\sum_{j=m_{t}+1}^{m_{t+1}}(\sum_{i=n_{j}+1}^{
n_{j+1}}  
c_{i}y_{i})]\|
$$\ 
$$
\leq {\frac{1}{K_t}}\|{\sum_{i=s}^{n_{l+1}}c_{i}y_{i}}\|
+ {\frac{1}{K_k}}\sum_{j=l+1}^{m_{k+1}}\|\sum_{i=n_{j}+1}^{n_{j+1}}  
c_{i}y_{i}\|
+ \sum_{t=k+1}^{\infty}\frac{1}{K_t}[\sum_{j=m_{t}+1}^{m_{t+1}}\|  
\sum_{i=n_{j}+1}^{n_{j+1}} c_{i}y_{i}]\|
$$ \

$$
\leq {\frac{2}{K_t}{\mbox{sup}}_{{\epsilon}_{i} = \pm  
1}\|\sum_{i=n_{l}+1}^{n_{l+1}} {\epsilon}_{i}c_{i}y_{i}\|} +
\sum_{j=l+1}^{\infty} \|\sum_{i=n_{j}+1}^{n_{j+1}} c_{i}y_{i}\| 
$$ \ $$ \leq {\frac{2}{K_t}} + \sum_{i=l}^{\infty} {\frac{1}{2^i}} \leq  
  \frac{2}{k} + \frac{1}{2^{l-1}}. $$
Hence, our series (8) converges. This contradiction
completes the proof of Theorem 3.2. \\ 
\begin{large}
\begin{center}
References:
\end{center}
\end{large}
\lbrack CC\rbrack \ Casazza, P.G. and Christensen, O.: {\it Frames
containing a Riesz basis and preservation of this property under 
perturbation.} Submitted, august 1995. \\ \\
\lbrack C1\rbrack \ Christensen, O.: {\it Frames containing a Riesz basis
and approximation of the frame coefficients using finite dimensional
methods.} Accepted for publication by J. Math. Anal. Appl.  \\ \\
\lbrack C2\rbrack \ Christensen, O.: {\it Frames and pseudo-inverse operators.}
J. Math. Anal. Appl. {\bf 195} (1995). \\ \\
\lbrack C3\rbrack \ Christensen, O.: {\it Frame perturbations.} Proc. Amer.
Math. Soc. {\bf 123} (1995) p.1217-1220. \\ \\
\lbrack CH\rbrack \ Christensen, O. and Heil, C.: {\it Perturbation of
Banach frames and atomic decomposition.} Accepted for publication by
Math. Nach..  \\ \\
\lbrack D\rbrack \ Diestel, J.: {\it Sequences and series in Banach spaces.}
Springer-Verlag, Graduate texts in mathematics, no. 92, 1984.  \\ \\
\lbrack H\rbrack \ Holub, J.: {\it Pre-frame operators, Besselian frames and
near Riesz bases.} Proc. Amer. Math. Soc. {\bf 122} (1994), p. 779-785. \\ \\
\lbrack He\rbrack \ Heil, C.: {\it Wavelets and frames.} Signal Processing,
Part 1, IMA Vol. Math. Appl. vol. 22. Springer-Verlag (1990) p.147-160. \\ \\
\lbrack Se\rbrack \ Seip, K.: {\it On the connection between exponential
bases and certain related sequences in $L^2(- \pi , \pi ) $.}
J. Funct. Anal. {\bf 130} (1995), p.131-160. 
\\ \\
{\bf Peter G. Casazza \\ Department of Mathematics \\ University of
Missouri \\ Columbia, Mo 65211 \\ USA \\ Email: pete@casazza.cs.missouri.edu
\\ \\ Ole Christensen \\ Mathematical Institute \\ Building 303 \\ 
Technical University of Denmark \\ 2800 Lyngby \\ Denmark \\
Email: olechr@mat.dtu.dk }

\end{document}